\documentclass[12pt,reqno]{amsart}
\usepackage{fullpage}
\usepackage{times}
\usepackage[utf8]{inputenc}
\usepackage[english]{babel}
\usepackage[mathscr]{euscript}
\usepackage{amsmath,amssymb,amsthm, graphicx, enumitem, xcolor}
\usepackage[colorlinks=true,citecolor=blue,linkcolor=blue,urlcolor=blue]{hyperref}
\newtheorem{theorem}{Theorem}[section]
\newtheorem{proposition}[theorem]{Proposition}
\newtheorem{lemma}[theorem]{Lemma}
\newtheorem{corollary}[theorem]{Corollary}
\newtheorem{question}[theorem]{Question}
\newtheorem{conjecture}[theorem]{Conjecture}
\newtheorem{claim}[theorem]{Claim}

\theoremstyle{definition}
\newtheorem{definition}[theorem]{Definition}

\newtheorem{remark}[theorem]{Remark}

\usepackage[capitalise]{cleveref}

\crefname{theorem}{Theorem}{Theorems}
\crefname{proposition}{Proposition}{Propositions}
\crefname{lemma}{Lemma}{Lemmas}
\crefname{corollary}{Corollary}{Corollaries}
\crefname{question}{Question}{Questions}
\crefname{conjecture}{Conjecture}{Conjectures}
\crefname{remark}{Remark}{Remarks}
\newenvironment{poc}{\begin{proof}[Proof of claim]}{\end{proof}}

\newcommand{\R}{\mathbb{R}} 
\newcommand{\N}{\mathbb{N}}

\title{The Erd\H{o}s similarity problem for polynomials and analytic maps}

\author{Alex Burgin}
\address{School of Mathematics\\ Georgia Institute of Technology\\ Atlanta\\ United States}
\email{alexander.burgin@gatech.edu}
\author{Chi Hoi Yip}
\address{Department of Mathematics, Hong Kong University of Science and Technology, Clear Water Bay, Hong Kong}
\email{machyip@ust.hk}


\begin{document}

\begin{abstract}
In this paper, we study the Erd\H{o}s similarity problem by replacing affine copies with images of polynomial maps or analytic maps. This is inspired by a connection between a partial converse of Szemer\'{e}di’s theorem and the Erd\H{o}s similarity problem and its variants, as well as a version of polynomial Szemer\'{e}di’s theorem by Bergelson and Leibman. In particular, we strengthen several known results in the literature. We show that a large class of sequence are not polynomially universal in this sense. First, for sequences satisfying suitable gap hypotheses, an explicit Cantor-type witness is constructed. Second, for sets satisfying certain growth conditions, a random Cantor-type witness is constructed, extending work of Kolountzakis. This also extends certain results in recent work of Cuéllar-Iosevich-Kulkarni-Aravena-Yavicoli to the polynomial setting. Finally, we show for arbitrary sequences $A$ with zero as an accumulation point that the $(d+3)$-fold sumset $(d+3)A$ is not polynomially universal when restricted to polynomials of degree at most $d$; this complements classical work of Bourgain.
\end{abstract}
\maketitle

\section{Introduction}
In 1920, Steinhaus \cite{S20} showed that if \(A\) is a finite subset of \(\R\), then every
measurable set \(E\subset\R\) of positive Lebesgue measure contains a nontrivial
affine copy of \(A\).  In 1974, Erd\H{o}s \cite{Erdos74} conjectured that this
finite-pattern phenomenon breaks down completely for infinite sets.  This is known as the \emph{Erd\H{o}s similarity conjecture}; see also Erd\H{o}s Problem~\#120
\cite{ErdosProblems120}.

\begin{conjecture}[Erd\H{o}s similarity conjecture]\label{qq}
If \(A\) is an infinite subset of \(\R\), then there is a set \(E\subset[0,1]\)
with \(|E|>0\) such that \(E\) does not contain any nontrivial affine copy of
\(A\).
\end{conjecture}

To study Conjecture~\ref{qq}, it suffices to consider sequences of positive real
numbers decreasing to \(0\) \cite{S00}.  Although the conjecture remains widely
open, several classes of sequences are known to satisfy the conjecture. 
Eigen~\cite{Eigen85} and Falconer~\cite{Falconer84} independently proved the conjecture
for sublacunary sequences, that is, sequences \(A=(a_n)_{n=1}^{\infty}\)
satisfying \(a_{n+1}/a_n\to1\). Humke and Laczkovich \cite{HumLac98} showed that the condition $\lim_{n \to \infty} a_{n+1}/a_n=1$ can be replaced with $\limsup_{n \to \infty} a_{n+1}/a_n=1$, and $a_m-a_{m+1} \leq a_n-a_{n+1}$ for all $m>n$. Kolountzakis
\cite{Kolountzakis97} proved a more flexible criterion in terms of large finite
subsets with large normalized gaps.  Bourgain \cite{Bourgain87} proved the
conjecture for sets with special additive structure, for instance sums of three
infinite sets.  We refer to \cite{S00,S24} for surveys of the Erd\H{o}s
similarity problem. Very recently, OpenAI~\cite{OpenAI26} claimed a proof of the conjecture for geometric progressions. See ~\cite{IY26, JX26} for expositions and extensions. Nevertheless, The Erd\H{o}s similarity conjecture remains open in general.

We study a variant of the Erd\H{o}s similarity problem in which affine copies are
replaced by polynomial copies, and more generally by analytic copies.  We now
explain the motivation for this variant.

The Erd\H{o}s similarity problem is naturally connected with finite-pattern
recurrence phenomena. As Bourgain~\cite{Bourgain87} observed, Erd\H{o}s' question was likely motivated by the finite-pattern phenomena underlying Szemer\'edi's theorem~\cite{Szemeredi75} and its extension by Furstenberg and Katznelson~\cite{FK78}. Szemer\'edi's theorem guarantees arbitrarily long arithmetic progressions in every set of integers of positive upper density, while Steinhaus' theorem implies that every measurable subset of $\mathbb R$ of positive measure contains a nontrivial affine copy of every finite configuration. The Erd\H{o}s similarity problem asks whether this universality breaks down for infinite configurations. The Erd\H{o}s similarity problem in the large is particularly reminiscent of Szemer'edi's theorem, as it asks whether a uniform largeness condition forces the presence of infinite patterns~\cite{BKM23,BGKMW23,KolountzakisPapageorgiou23, GMY}. It was recently resolved by Yip~\cite{Y26}, who proved that for every unbounded set $A\subset\mathbb R$ and every $0<\varepsilon<1$, there exists a closed set $S\subset\mathbb R$ having measure at least $1-\varepsilon$ in every unit interval and containing no nontrivial affine copy of $A$.

The polynomial version of the problem is motivated by the corresponding polynomial
recurrence theory.  A real-line polynomial Szemer\'edi theorem of Bergelson and
Leibman \cite[Theorem B'']{BL96} says that if \(P\in\R[x]\) satisfies \(P(0)=0\),
if \(a_1,\ldots,a_k\in\R\), and if \(E\subset\R\) has positive upper Banach
density, then there are a positive integer \(n\) and a point \(x\in\R\) such that
\[
        x+P(na_j)\in E\qquad(1\leq j\leq k).
\]
Equivalently, with
\[
        Q(t)=x+P(nt),
\]
the set \(E\) contains the finite polynomial image
\[
        \{Q(a_1),\ldots,Q(a_k)\}.
\]
Thus finite polynomial copies of arbitrary finite configurations are unavoidable
in positive-density sets.  We ask whether infinite polynomial copies can
nevertheless be avoided by a set of positive measure.

This perspective is part of a broader circle of continuous analogues of classical
results from additive combinatorics. For example, polynomial Roth-type theorems on sets of
fractional dimension have been studied in \cite{DGR19,FGP22}, and continuous
Furstenberg--S\'{a}rk\"{o}zy-type problems have been studied in
\cite{DKS23,KOS23}.  A more familiar polynomial Szemer\'edi theorem of Bergelson
and Leibman \cite[Theorem B']{BL96} concerns configurations \(u+p_i(n)\) generated
by several polynomials \(p_1,\ldots,p_k\).  By the Lebesgue density theorem, this
also has a continuous Steinhaus-type analogue; see for instance
\cite[p.~6973]{DGR19}.

It is therefore natural to ask for an Erd\H{o}s similarity analogue in which
affine maps are replaced by polynomial maps, and more generally by analytic maps.
This direction is also motivated by the bi-Lipschitz variant studied recently by
Feng, Lai, and Xiong \cite{FLX23}: since every nontrivial affine map is
bi-Lipschitz, their problem already contains the classical affine setting and asks
for avoidance under a larger class of transformations.

More precisely, we study the following polynomial and analytic analogue of the
Erd\H{o}s similarity problem.

\begin{question}[Polynomial and analytic Erd\H{o}s similarity problem]\label{polysmall}
Let \((a_n)_{n=1}^{\infty}\) be a decreasing sequence of positive numbers with
\(a_n\to0\).  Does there exist a measurable set \(E\subset[0,1]\) with
\(|E|>0\) such that
\[
        \{P(a_n):n\in\N\}\not\subset E
\]
for every nonconstant polynomial \(P\in\R[x]\)?  Can one choose \(E\) to avoid
the images of \((a_n)\) under every nonconstant real-analytic function defined
near \(0\)?
\end{question}

We prove positive results for Question~\ref{polysmall} by three complementary methods.  The first is an \emph{explicit Cantor-type construction}.  It is deterministic, works particularly well for sequences satisfying a global gap-comparability condition, and gives a strong tail-avoidance conclusion with a set of measure arbitrarily close to one.  

We first state the result obtained by the explicit construction.
\begin{theorem}\label{thm:analytic}
Let \((a_n)_{n=1}^{\infty}\) be a decreasing sequence of positive numbers such that \(a_n\to0\),
\[
        \limsup_{n\to\infty}\frac{a_{n+1}}{a_n}=1,
\]
and there is a constant \(C>0\) such that
\[
        a_m-a_{m+1}\le C(a_n-a_{n+1})\qquad(m>n).
\]
Then for each \(\epsilon>0\), there is a set \(E\subset[0,1]\) with \(|E|\ge1-\epsilon\) such that \(E\) does not contain
\[
        \{f(a_n):n\ge N\}
\]
for any nonconstant function \(f:\R\to\R\) analytic at \(0\) and any positive integer \(N\).
\end{theorem}

The second is a \emph{random Cantor-type construction}, in the spirit of Kolountzakis \cite{Kolountzakis97}.  It is less explicit but more flexible: it applies under a local finite-block hypothesis, treats polynomial maps by counting finitely many interval patterns, and treats analytic maps by adding a Taylor-truncation estimate. The following theorem is a strengthening of a result of Kolountzakis \cite[Theorem 3]{Kolountzakis97}.

\begin{theorem}\label{thm:randommain}
Let
\[
        A=(a_n)_{n=1}^{\infty},\qquad
        a_1>a_2>\cdots>0,\qquad a_n\to0,
\]
and identify \(A\) with its underlying set. Suppose that there are finite subsets
\[
        B_\nu=\{b_{\nu,1}>\cdots>b_{\nu,N_\nu}>0\}\subset A
\]
such that
\[
        N_\nu\to\infty,\qquad b_{\nu,1}\to0,
\]
and, writing
\[
        \Delta_\nu
        =
        \min_{1\le i<N_\nu}
        \frac{b_{\nu,i}-b_{\nu,i+1}}{b_{\nu,1}},
\]
one has
\[
        -\log\Delta_\nu=o(N_\nu).
\]
Then there exists a measurable set \(E\subset[0,1]\) of positive Lebesgue measure
such that
\[
        P(A)\not\subset E
\]
for every nonconstant real polynomial \(P\).

Moreover, if the sets \(B_\nu\) may be chosen so that
\begin{equation}\label{eq:entropy-main}
        \left(
        1+\frac{|\log\Delta_\nu|}{\log(1/b_{\nu,1})}
        \right)
        \bigl(\log N_\nu+|\log\Delta_\nu|\bigr)
        =
        o(N_\nu),
\end{equation}
then \(E\) may be chosen so that
\[
        f(A)\not\subset E
\]
for every nonconstant real-analytic function \(f\) defined in a neighbourhood of
\(0\).
\end{theorem}

\begin{remark}
This formulation is equivalent, for the polynomial conclusion, to Kolountzakis's
large-gap hypothesis that \(A\) contains arbitrarily large finite subsets \(B\) with
\(-\log\Delta(B)=o(|B|)\). Indeed, since \(A\) is a sequence tending to \(0\), one
may pass to terminal portions of such good blocks to ensure that their largest
elements tend to \(0\), without worsening the normalized gap condition.
\end{remark}

The following consequence gives convenient hypotheses which imply the analytic part of the random Cantor theorem.  In particular, it contains the sublacunary and gap-comparable classes mentioned above.

\begin{corollary}\label{cor:sufficient-criteria}
Let
\[
        A=(a_n)_{n=1}^\infty,
        \qquad a_1>a_2>\cdots>0,
        \qquad a_n\to0,
\]
and identify the sequence with its underlying set.  Assume that at least one of the following conditions holds.
\begin{enumerate}
\item[\rm(i)] For each \(n\), let \(\ell(n)\) be the smallest integer such that
\[
        a_{\ell(n)}\leq \frac{a_n}{4},
\]
and put \(T(n)=\{n,n+1,\ldots,\ell(n)\}\). Then
\[
        \liminf_{n\to\infty}
        \frac{\max_{m\in T(n)}(a_m-a_{m+1})}{a_n}=0.
\]

\item[\rm(ii)]
\[
        \frac{a_{n+1}}{a_n}\to1.
\]

\item[\rm(iii)]
\[
        \limsup_{n\to\infty}\frac{a_{n+1}}{a_n}=1,
\]
and there exists \(C>0\) such that
\[
        a_m-a_{m+1}\leq C(a_n-a_{n+1})
        \qquad(m>n).
\]
\end{enumerate}
Then there exists a measurable set \(E\subset[0,1]\) of positive Lebesgue measure
such that
\[
        f(A)\not\subset E
\]
for every nonconstant real-analytic function \(f\) defined in a neighbourhood of
\(0\).
\end{corollary}

We now discuss another corollary to this theorem; this crucially uses an estimate in \cite{cuellar2026nonuniversalitysumsetslacunarysequences}, which shows that certain lacunary sequences cannot have too large additive energy.

Let us begin by setting up the framework. We call a sequence $(b_i)_{i\geq 1}\subset (0,\infty)$ \textit{lacunary with ratio} $q\in (0,1)$ if $b_{i+1}\leq qb_i$ for all $i$, and its associated counting function is \begin{align*}
    I(W):=\#\{i:b_i\geq \exp(-W)\},\quad W>0.
\end{align*} 

\begin{corollary}\label{cor:lacunarySumsets}
Suppose $S_1$ and $S_2$ are two lacunary sequences with counting functions $I(W)$ and $J(W)$, respectively. Then if \begin{align*} 
\limsup_{W\rightarrow \infty}\frac{I(W)J(W)}{W}=\infty
\end{align*} then neither $S_1+S_2$ nor $S_1-S_2$ are polynomial universal.
\end{corollary}

Finally, restricting to polynomials of fixed degree, we are able to show a general result for sufficiently large sumsets of an arbitrary infinite sequence. For $r\geq 1$, let $rA:=\{a_1+\cdots+a_r:a_i\in A\}$. We can show the following result.

\begin{theorem}\label{thm:sumsetNonuniversality}
    Let $d\geq 1$ be an integer and suppose $A\subset (0,\infty)$ has zero as an accumulation point. For every $\epsilon\in (0,1)$, there is a Borel set $E\subset [0,1]$ with Lebesgue measure $\mu(E)>1-\epsilon$ and \begin{align*}
        P((d+3)A)\not\subset E\quad \text{for every nonconstant }P\in \mathbb{R}[y]\text{ with }\deg P\leq d.
    \end{align*}
\end{theorem}

Our proof in \S\ref{sec:fixedDegreeNU} resembles Bourgain's proof \cite{Bourgain87} for the $d=1$ case and Tao's related exposition \cite{Tao2021Toolkit}, but extends this to degree-$d$ polynomial non-universality.

\section{The explicit Cantor construction}

This section proves Theorem~\ref{thm:analytic}.  The construction is deterministic: at each scale we remove regularly spaced gaps and then intersect the resulting high-measure sets.  The proof is inspired by a recent paper by Feng, Lai, and Xiong \cite{FLX23}, where they studied the Erd\H{o}s similarity problem via bi-Lipschitz embeddings. 

For our purpose, we need the following strengthening of \cite[Theorem 1.3]{FLX23}. We say a function $f: \R \to \R$ is \emph{bi-Lipschitz near $0$} if there exists $\delta>0$ and $L>1$ such that
$$
L^{-1}|x-y|\leq |f(x)-f(y)|\leq L|x-y|, \quad \forall x,y \in [-\delta, \delta].
$$
\begin{lemma}\label{lem:bi-Lispschitz}
Let $(a_n)_{n=1}^{\infty}$ be a decreasing sequence of positive numbers such that $\lim_{n \to \infty} a_n=0$, $\limsup_{n \to \infty} a_{n+1}/a_n=1$, and there is a constant $C>0$ such that $a_m-a_{m+1} \leq C(a_n-a_{n+1})$ for all $m>n$. Then for each $\epsilon \in (0,1)$, there is a set $E \subset [0,1]$ with $|E| \geq 1-\epsilon$, such that $E$ does not contain $\{f(a_n): n \geq N\}$ for any function \(f\) defined in a neighbourhood of \(0\) and bi-Lipschitz near \(0\), and any positive integer \(N\).       
\end{lemma}
\begin{proof}
Since $\limsup_{n \to \infty} a_{n+1}/a_n=1$, we have $\liminf_{n \to \infty}(a_{n}-a_{n+1})/a_n=0$. In particular, we can find a strictly increasing sequence $(n_k)_{k=1}^{\infty}$ such that
\begin{equation}\label{n_k}
\frac{a_{n_k}-a_{n_k+1}}{a_{n_k}}<\frac{\epsilon}{2^{k}k^2}.    
\end{equation}
Let $\delta_k=k(a_{n_k}-a_{n_{k}+1})$. Let $\ell_k$ be the smallest positive integer such that
$$
\frac{2a_{n_k}}{k}\geq \frac{1}{\ell_k};
$$
then we have $1/\ell_k \geq a_{n_k}/k$. 
It then follows from inequality~\eqref{n_k} that
$$
\delta_k \ell_k\leq  k^2\frac{a_{n_k}-a_{n_{k}+1}}{a_{n_k}} <\frac{\epsilon}{2^k}.
$$
Let 
$$E_k=\bigcup_{j=0}^{\ell_k-1}\bigg[\frac{j}{\ell_k}, \frac{j+1}{\ell_k}-\delta_k\bigg].$$ 
Then $|E_k|=1-\delta_k \ell_k\geq 1-\epsilon/2^k$. Let $E=\cap_{k=1}^{\infty} E_k$. Then we have $|E|\geq 1-\epsilon$. 

We claim that $E$ has the desired property. For the sake of contradiction, assume that 
$f: \R \to \R$ is a function bi-Lipschitz near $0$ and $N$ is a positive integer, such that $\{f(a_n): n \geq N\} \subset E$. Then there exists $\delta>0$ and $L>1$ such that
\begin{equation}\label{biLip}
L^{-1}|x-y|\leq |f(x)-f(y)|\leq L|x-y|, \quad \forall x,y \in [-\delta, \delta].    
\end{equation}   
By possibly increasing $N$, we may assume that $a_N \leq \delta$. Note that $f$ is monotone in $[-\delta, \delta]$.  Let $k$ be an integer such that $k>2CL$ and $n_k>N$. Then inequality~\eqref{biLip} implies that
$$
|f(a_{n_k})-f(0)|=\lim_{m \to \infty} |f(a_{n_k})-f(a_m)| \geq \lim_{m \to \infty} L^{-1}(a_{n_k}-a_m)=\frac{a_{n_k}}{L}>\frac{2a_{n_k}}{k}\geq \frac{1}{\ell_k}.
$$
Thus, $B=\{f(a_n): n \geq n_k\}$ is a subset of $E_k$ with diameter at least $\frac{1}{\ell_k}$. Recall that $E_k$ is a disjoint union of intervals of length less than $\frac{1}{\ell_k}$, so $B$ intersects at least two of those intervals forming $E_k$. On the other hand, for each $m \geq n_k$, inequality~\eqref{biLip} implies that
$$
|f(a_{m})-f(a_{m+1})|\leq L(a_{m}-a_{m+1}) \leq CL(a_{n_k}-a_{n_k+1})<k(a_{n_k}-a_{n_k+1})=\delta_k.
$$
However, the gap between any two intervals forming $E_k$ has size at least $\delta_k$, a contradiction. This finishes the proof of the claim.
\end{proof}

Now we are ready to prove Theorem~\ref{thm:analytic}.
\begin{proof}[Proof of Theorem~\ref{thm:analytic}]
Let $\epsilon>0$. For each positive integer $k$, we have $\limsup_{n \to \infty} a_{n+1}^k/a_n^k=1$; moreover, for all $m>n$, we have
\begin{align*}
a_m^k-a_{m+1}^k
&=(a_m-a_{m+1})\bigg(\sum_{j=0}^{k-1}a_m^ja_{m+1}^{k-1-j}\bigg) \leq ka_m^{k-1}(a_m-a_{m+1}) \leq ka_{n+1}^{k-1}(a_m-a_{m+1})\\
&\leq Cka_{n+1}^{k-1}(a_n-a_{n+1}) \leq C(a_n-a_{n+1})\bigg(\sum_{j=0}^{k-1}a_n^ja_{n+1}^{k-1-j}\bigg)\leq C(a_n^k-a_{n+1}^k).
\end{align*}
Thus we can apply Lemma~\ref{lem:bi-Lispschitz} to the sequence $(a_n^k)_{n=1}^{\infty}$ to obtain a set $E_k \subset [0,1]$ with $|E_k|\geq 1-\epsilon/2^k$, such that $E_k$ does not contain $\{f(a_n^k): n \geq N\}$ for any  function $f: \R \to \R$ bi-Lipschitz near $0$
and any positive integer $N$. 

Let $E=\cap_{k=1}^{\infty} E_k$. Then we have $|E|\geq 1-\epsilon$. 
We claim that $E$ has the desired property. Suppose otherwise that there is a non-constant function $f: \R \to \R$ analytic at 0 and a positive integer $N$, such that $\{f(a_n): n \geq N\} \subset E$. Since $f$ is analytic at $0$, there is $\delta \in (0,1)$, such that we can write 
$$f(x)=\sum_{i=0}^{\infty} c_ix^i$$
and the series converges to $f(x)$ for $|x|\leq \delta$. Since $f$ is non-constant, we can find the smallest positive integer $k$ such that $c_k \neq 0$.
Since $a_n \to 0$ as $n \to \infty$, there is $m \geq N$ such that $a_n \leq \delta$ for all $n \geq m$. For $0\leq y\leq \delta^k$, put
$$
        g_+(y)=f(y^{1/k}).
$$
Since $k$ is the first positive integer with $c_k\neq0$, we have
$$
        g_+(y)=c_0+c_ky+O(y^{1+1/k})
        \qquad (y\downarrow0).
$$
Define $g$ on a neighbourhood of $0$ by
$$
        g(y)=g_+(y)\quad (y\geq0),
        \qquad
        g(y)=2c_0-g_+(-y)\quad (y<0).
$$
Then
$$
        g(y)=c_0+c_ky+O(|y|^{1+1/k})
        \qquad (y\to0),
$$
and hence $g$ is bi-Lipschitz near $0$. On the other hand,
$f(a_n)=g(a_n^k)$ for all $n\geq m$. It follows that
$$
\{g(a_n^k): n \geq m\}=\{f(a_n): n \geq m\} \subset \{f(a_n): n \geq N\} \subset E \subset E_k,
$$
which violates the defining property of $E_k$. This finishes the proof of the claim.
\end{proof}

\section{The random Cantor construction}\label{sec:random}
The second method is a random Cantor-type construction.  It is less explicit than the construction in the previous section, but it is more flexible: it only requires large finite blocks with good normalized gaps.  The main point is that, after a random partition of \([0,1]\), the uncountable family of polynomial or analytic maps gives only subexponentially many interval patterns on the chosen finite block.

\subsection{A lemma on hyperplane intersection}

We prove the following elementary counting fact.

\begin{lemma}\label{lem:threshold-count}
Let $H\geq 1$. Let \(L_1,\ldots,L_H\) be affine-linear functions on \(\R^r\).  Then the number
of possible truth-value patterns
\[
        \bigl(\mathbf 1_{L_1(c)\ge0},\ldots,\mathbf 1_{L_H(c)\ge0}\bigr),
        \qquad c\in\R^r,
\]
is at most \((r+1)H^r\).
\end{lemma}

\begin{proof}
Let \(N_r(H)\) denote the maximum possible number of such truth-value patterns
cut out by \(H\) affine inequalities in \(\R^r\).  We prove by induction that
\[
        N_r(H)\le \sum_{q=0}^r \binom{H}{q}.
\]

Consider adding the last inequality \(L_H(c)\ge0\).  If \(L_H\) is constant, then
it does not increase the number of patterns.  Otherwise, its zero set
\[
        Z=\{c\in\R^r:L_H(c)=0\}
\]
is an affine subspace of dimension \(r-1\).  The only patterns for
\(L_1,\ldots,L_{H-1}\) that can split into two different patterns after adding
\(L_H\) are those whose regions meet \(Z\).  The number of such possible boundary
patterns is at most the number of truth-value patterns cut out on \(Z\) by the
restricted affine functions
\[
        L_1|_Z,\ldots,L_{H-1}|_Z.
\]
Since \(Z\) has dimension \(r-1\), this number is at most \(N_{r-1}(H-1)\).  Hence
\[
        N_r(H)\le N_r(H-1)+N_{r-1}(H-1).
\]
Together with \(N_0(H)\le1\) and \(N_r(0)=1\), induction gives
\[
        N_r(H)\le \sum_{q=0}^r \binom{H}{q}\leq (r+1)H^r. \qedhere
\]
\end{proof}

\subsection{Proof of the polynomial assertion}

\begin{lemma}\label{lem:polynomial-window}
Fix integers \(D\geq m\geq1\) and real numbers \(0<b<B<\infty\). Let
\[
        \mathcal W_{D,m,b,B}
        =
        \left\{
        c=(c_m,\ldots,c_D):
        b\leq |c_m|\leq B,\quad |c_j|\leq B\text{ for }j>m
        \right\}.
\]
For \(c\in\mathcal W_{D,m,b,B}\), set
\[
        R_c(t)=c_mt^m+c_{m+1}t^{m+1}+\cdots+c_Dt^D.
\]
Let
\[
        A_\nu=\{a_{\nu,1}>\cdots>a_{\nu,n_\nu}>0\}
\]
be finite subsets satisfying
\[
        a_{\nu,1}\to0,
        \qquad
        -\log\delta_\nu=o(n_\nu),
        \qquad
        \delta_\nu=
        \min_{1\leq i<n_\nu}
        \frac{a_{\nu,i}-a_{\nu,i+1}}{a_{\nu,1}}.
\]
Then, for every \(\varepsilon,\eta \in (0,1)\), and for all sufficiently large \(\nu\), there
is a measurable set \(F\subset[0,1]\) with
\[
        \mu(F)>1-\varepsilon
\]
such that the exceptional set
\[
        X_\nu(F)
        =
        \left\{
        \gamma\in[0,1]:
        \exists c\in\mathcal W_{D,m,b,B}
        \text{ such that }
        \gamma+R_c(a_{\nu,i})\in F
        \text{ for all }1\leq i\leq n_\nu
        \right\}
\]
satisfies
\[
        \mu(X_\nu(F))<\eta.
\]
\end{lemma}

\begin{proof}
We first record a uniform
monotonicity estimate for $R_c$ near \(0\), uniform for $c\in \mathcal W_{D,m,b,B}$. Choose \(\rho>0\), depending only on \(D,m,b,B\), so small that
\[
        B\sum_{j=m+1}^D j\rho^{j-m}\leq \frac{mb}{2}.
\]
\begin{claim}\label{claim:sep}
For each $c\in \mathcal W_{D,m,b,B}$, \(R_c\) is monotone on \((0,\rho]\), and for \(0<y<x\leq\rho\),
\[
        |R_c(x)-R_c(y)| \geq \frac b2(x-y)^m.
\]
\end{claim}
\begin{poc}
Fix $c\in \mathcal W_{D,m,b,B}$. Then, for \(0<t\leq\rho\), we have
\[
        R_c'(t)=mc_mt^{m-1}+\sum_{j=m+1}^D jc_jt^{j-1},
\]
where
\[
\left|
\sum_{j=m+1}^D jc_jt^{j-1}
\right|
\le
t^{m-1}\sum_{j=m+1}^D j|c_j|t^{j-m}  
\le
B t^{m-1}\sum_{j=m+1}^D j\rho^{j-m}
\le
\frac{mb}{2}t^{m-1}.
\]
Thus, for \(0<t\leq\rho\),
\[
        |R_c'(t)|\geq \frac{mb}{2}t^{m-1},
\]
and \(R_c'\) has constant sign. Hence \(R_c\) is monotone on \((0,\rho]\), and for \(0<y<x\leq\rho\),
\[
        |R_c(x)-R_c(y)|
        \geq
        \frac{mb}{2}\int_y^x t^{m-1}\,dt
        =
        \frac b2(x^m-y^m)
        \geq
        \frac b2(x-y)^m. \qedhere
\]
\end{poc}
Fix a large \(\nu\), and write
\[
        n=n_\nu,\qquad a_i=a_{\nu,i},\qquad \delta=\delta_\nu.
\]
We may assume \(a_1<\rho\). Let
\[
        r=D-m+1,
        \qquad
        k=
        \left\lceil
        \frac{8}{b(a_1\delta)^m}
        \right\rceil.
\]
Partition \([0,1)\) into the half-open intervals
\[
        I_\ell=\left[\frac{\ell}{k},\frac{\ell+1}{k}\right)
        \quad(0\leq \ell\leq k-1).
\]
We keep each interval independently with probability \(p=p_\nu\) to be determined later in the proof, and let \(F\) be
the resulting random set.  For each \(x\in[0,1)\), write \(I(x)\) for the unique interval
of this partition containing \(x\).

 Fix
\(\gamma\in[0,1]\).  For \(c\in\mathcal W_{D,m,b,B}\), suppose first that
\[
        \gamma+R_c(a_i)\in[0,1)\qquad(1\le i\le n).
\]
Then define the interval pattern of \(c\) by
\[
        \pi_\gamma(c)=(\ell_1(c),\ldots,\ell_n(c)),
\]
where
\[
        \gamma+R_c(a_i)\in I_{\ell_i(c)}.
\]
Since \(R_c\) is monotone on \((0,\rho]\), the values
\(R_c(a_1),\ldots,R_c(a_n)\) are ordered.  For each \(1\leq i<n\), Claim~\ref{claim:sep} implies that
\[
        |R_c(a_i)-R_c(a_{i+1})|
        \geq
        \frac b2(a_i-a_{i+1})^m
        \geq
        \frac b2(a_1\delta)^m
        >
        \frac1k.
\]
Thus the intervals
\[
        I_{\ell_1(c)},\ldots,I_{\ell_n(c)}
\]
are distinct.

We now count the number of possible interval patterns as \(c\) varies in \(\mathcal W_{D,m,b,B}\). 
We claim that only
\[
        O_{b,B,D}\bigl(n^r\delta^{-mr}\bigr)
\]
different interval patterns can occur.  For fixed \(i\), the interval containing
\(\gamma+R_c(a_i)\) is determined by which partition endpoints lie to its left.
Equivalently, it is determined by the truth values of the affine inequalities
\[
        \gamma+R_c(a_i)\ge \frac{\ell}{k}.
\]
Writing this out, these inequalities are
\[
        \gamma+c_ma_i^m+c_{m+1}a_i^{m+1}+\cdots+c_Da_i^D
        \ge
        \frac{\ell}{k}.
\]

For fixed \(i\), the set of possible values of \(R_c(a_i)\), with
\(c\in\mathcal W_{D,m,b,B}\), is contained in an interval of length
\(O_{B,D}(a_i^m)\), hence in an interval of length \(O_{B,D}(a_1^m)\). Since the partition
endpoints are spaced by \(1/k\), only
\[
        O_{b,B,D}(ka_1^m+1)
        =
        O_{b,B,D}(\delta^{-m})
\]
values of \(\ell\) can change their truth value as \(c\) varies. Thus, as \(i\) ranges from \(1\) to \(n\), the interval pattern is determined by
the truth values of at most
\[
        H\ll_{b,B,D}n\delta^{-m}
\]
affine inequalities in the \(r=D-m+1\) coefficient variables.  By
Lemma~\ref{lem:threshold-count}, the number of possible truth-value patterns,
and hence the number of possible interval patterns \(\pi_\gamma(c)\), is at most
\[
        (r+1)H^r
        \ll_{b,B,D}
        n^r\delta^{-mr}.
\]

Now, if \(\gamma\in X_\nu(F)\), then one of these possible interval patterns has
all of its \(n\) intervals kept.  Since the intervals in each pattern are distinct,
the probability that a fixed pattern survives is \(p^n\).  Therefore, uniformly in
\(\gamma\),
\[
        \mathbb P(\gamma\in X_\nu(F))
        \le
        C_{b,B,D}n^r\delta^{-mr}p^n.
\]
After integrating over \(\gamma\in[0,1]\), we get
\[
        \mathbb E\,\mu(X_\nu(F))
        \leq
        C_{b,B,D}n^r\delta^{-mr}p^n.
\]

Put
\[
        S_\nu=r\log n_\nu+mr|\log\delta_\nu|.
\]
Since \(m\) and \(r\) are fixed and \(-\log\delta_\nu=o(n_\nu)\), we have
\[
        S_\nu=o(n_\nu).
\]
Choose
\[
        p_\nu=\exp\left(-\sqrt{\frac{S_\nu+1}{n_\nu}}\right).
\]
Then \(p_\nu\to1\). Moreover,
\[
\begin{aligned}
\log\left(n_\nu^r\delta_\nu^{-mr}p_\nu^{n_\nu}\right)
&=
r\log n_\nu+mr|\log\delta_\nu|
-
n_\nu\sqrt{\frac{S_\nu+1}{n_\nu}}  \\
&=
S_\nu-\sqrt{n_\nu(S_\nu+1)}
\to -\infty,
\end{aligned}
\]
because \(S_\nu+1=o(n_\nu)\). Hence
\[
        n_\nu^r\delta_\nu^{-mr}p_\nu^{n_\nu}\to0.
\]
In particular, for sufficiently large \(\nu\), we have
\[
        1-p_\nu<\varepsilon/4,
        \qquad
        \mathbb E\,\mu(X_\nu(F))<\eta/4.
\]
Since \(\mathbb E\,\mu(F)=p_\nu\), Markov's inequality gives
\[
        \mathbb P(\mu(F)\leq1-\varepsilon)=\mathbb P(1-\mu(F)\geq\varepsilon)<\frac14,
        \qquad
        \mathbb P(\mu(X_\nu(F))\geq\eta)<\frac14.
\]
Thus, with positive probability, both desired inequalities hold. Fix such a
realization of \(F\). This proves the lemma.
\end{proof}

We are now ready to prove \cref{thm:randommain}.  

\begin{proof}[Proof of the polynomial assertion in \cref{thm:randommain}]
 Consider finite blocks
\[
        A_\nu=\{a_{\nu,1}>\cdots>a_{\nu,n_\nu}>0\}\subset A
\]
satisfying
\[
        a_{\nu,1}\to0,
        \qquad
        -\log\delta_\nu=o(n_\nu).
\]

Enumerate all quadruples
\[
        (D_j,m_j,b_j,B_j),
        \qquad
        D_j\geq1,
        \quad
        1\leq m_j\leq D_j,
        \quad
        0<b_j<B_j,
\]
where \(D_j,m_j\) are integers and \(b_j,B_j\in\mathbb Q\). This enumeration covers
every nonzero polynomial coefficient vector: if
\[
        P(t)=\gamma+c_1t+\cdots+c_Dt^D
\]
is nonconstant and \(D=\deg P\), then, with
\[
        m=\min\{i\geq1:c_i\neq0\},
\]
there is some \(j\) with \(D_j=D\), \(m_j=m\), and
\[
        (c_m,\ldots,c_D)\in\mathcal W_{D_j,m_j,b_j,B_j}.
\]
Choose positive numbers \(\varepsilon_j,\eta_j\) such that
\[
        \sum_j\varepsilon_j<\frac14,
        \qquad
        \sum_j\eta_j<\frac18.
\]
For each \(j\), apply Lemma~\ref{lem:polynomial-window} to
\(\mathcal W_{D_j,m_j,b_j,B_j}\), with parameters \(\varepsilon_j,\eta_j\). We obtain
a block \(A_{\nu_j}\subset A\), a measurable set \(F_j\subset[0,1]\), and an
exceptional set \(X_j\subset[0,1]\) such that
\[
        \mu(F_j)>1-\varepsilon_j,
        \qquad
        \mu(X_j)<\eta_j.
\]
Put
\[
        F=\bigcap_{j=1}^\infty F_j,
        \qquad
        X=\bigcup_{j=1}^\infty X_j.
\]
Then \(\mu(F)>3/4\) and \(\mu(X)<1/8\). Choose an open set \(G\supset X\) with
\(\mu(G)<1/4\), and set
\[
        E=F\setminus G.
\]
Then \(E\subset[0,1]\) has positive measure.

Suppose that \(P(A)\subset E\) for some nonconstant real polynomial \(P\). Write
\[
        P(t)=\gamma+c_1t+\cdots+c_Dt^D,
        \qquad D=\deg P,
\]
and let \(m=\min\{i\geq1:c_i\neq0\}\). Since \(a_n\to0\) and
\(E\subset[0,1]\), we have \(P(a_n)\to P(0)=\gamma\), so \(\gamma\in[0,1]\). Choose
\(j\) with
\[
        D_j=D,
        \qquad
        m_j=m,
        \qquad
        (c_m,\ldots,c_D)\in\mathcal W_{D_j,m_j,b_j,B_j}.
\]
The block \(A_{\nu_j}\) used in constructing \(F_j\) is a subset of \(A\), and
\(P(A)\subset E\subset F_j\). Hence \(P(A_{\nu_j})\subset F_j\), so by the definition
of \(X_j\) we have \(\gamma\in X_j\subset G\). But \(G\) is open and
\(P(a_n)\to\gamma\) as $n\to \infty$. Thus \(P(a_n)\in G\) for all sufficiently large $n$, contradicting \(P(A)\subset E=F\setminus G\). 
\end{proof}
\subsection{Proof of the analytic assertion}

\begin{lemma}\label{lem:analytic-window}
Fix
\[
        m\geq1,
        \qquad
        R>0,
        \qquad
        M>0,
        \qquad
        0<b<B<\infty.
\]
Let \(\mathcal V_{m,R,M,b,B}\) be the family of functions
\[
        h(t)=c_mt^m+c_{m+1}t^{m+1}+\cdots
\]
with real coefficients, satisfying
\[
        b\le |c_m|\le B,\qquad
        \sum_{j=m}^{\infty}|c_j|R^j\le M.
\]
Let
\[
        A_\nu=\{a_{\nu,1}>\cdots>a_{\nu,n_\nu}>0\}
\]
be finite blocks such that
\[
        a_{\nu,1}\to0,
        \qquad
        -\log\delta_\nu=o(n_\nu),
        \qquad
        \delta_\nu=
        \min_{1\leq i<n_\nu}
        \frac{a_{\nu,i}-a_{\nu,i+1}}{a_{\nu,1}},
\]
and
\begin{equation}\label{eq:entropy-blocks}
\left(
        1+\frac{|\log\delta_\nu|}{\log(1/a_{\nu,1})}
        \right)
        \bigl(\log n_\nu+|\log\delta_\nu|\bigr)
        =o(n_\nu).    
\end{equation}
Then, for every \(\varepsilon,\eta \in (0,1)\), and for all sufficiently large \(\nu\), there
is a measurable set \(F\subset[0,1]\) with
\[
        \mu(F)>1-\varepsilon
\]
such that the exceptional set
\[
        X_\nu(F)=
        \left\{
        \gamma\in[0,1]:
        \exists h\in\mathcal V_{m,R,M,b,B}
        \text{ such that }
        \gamma+h(a_{\nu,i})\in F
        \text{ for all }1\le i\le n_\nu
        \right\}
\]
satisfies \(\mu(X_\nu(F))<\eta\).
\end{lemma}

\begin{proof}
We first record a uniform
monotonicity estimate for $h$ near \(0\), uniform for $h\in \mathcal V_{m,R,M,b,B}$. Choose \(0<\rho<R\), depending only on \(m,R,M,b\), so small that
\[
        R^{-m}M\sup_{j\ge m+1}j\left(\frac{\rho}{R}\right)^{j-m}
        \le \frac{mb}{2}.
\]

\begin{claim}\label{claim:analytic}
Every \(h\in\mathcal V_{m,R,M,b,B}\) is monotone on
\((0,\rho]\), and for \(0<y<x\le\rho\),
\[
        |h(x)-h(y)|
        \ge
        \frac b2(x-y)^m.
\]
\end{claim}

\begin{poc}
Write
\[
        h(t)=\sum_{j=m}^{\infty}c_jt^j\in\mathcal V_{m,R,M,b,B}.
\]
Since
\[
        \sum_{j=m}^{\infty}|c_j|R^j\le M,
\]
for \(0<t\le\rho\) we have
\[
\begin{aligned}
\left|
\sum_{j=m+1}^{\infty}jc_jt^{j-1}
\right|
\le
t^{m-1}R^{-m}
\sum_{j=m+1}^{\infty}j|c_j|R^j
\left(\frac{t}{R}\right)^{j-m}  \le
t^{m-1}R^{-m}M \sup_{j\geq m+1}
j \left(\frac{t}{R}\right)^{j-m}  
\le
\frac{mb}{2}t^{m-1}.
\end{aligned}
\]
Since \(|c_m|\ge b\), it follows that \(h'(t)\) has the sign of \(c_m\) and
\[
        |h'(t)|\ge \frac{mb}{2}t^{m-1}
        \qquad(0<t\le\rho).
\]
Integrating from \(y\) to \(x\) gives
\[
        |h(x)-h(y)|
        \ge
        \frac{mb}{2}\int_y^x t^{m-1}\,dt
        =
        \frac b2(x^m-y^m) \geq \frac b2(x-y)^m. \qedhere
\]
\end{poc}

Fix a large \(\nu\), and write
\[
        n=n_\nu,\qquad a_i=a_{\nu,i},\qquad \delta=\delta_\nu.
\]
We may assume
\[
        a_1<\rho,\qquad a_1<R/4.
\]
Let
\[
        k=
        \left\lceil
        \frac{16}{b(a_1\delta)^m}
        \right\rceil.
\]
Partition \([0,1)\) into the intervals
\[
        I_\ell=\left[\frac{\ell}{k},\frac{\ell+1}{k}\right),
        \qquad 0\leq \ell\leq k-1.
\]
Choose \(0<\theta<1/4\) so small that
\[
        2\theta<\varepsilon/4.
\]
For each \(\ell\), let
\[
        J_\ell=
        \left[\frac{\ell+\theta}{k},\frac{\ell+1-\theta}{k}\right)
        \subset I_\ell .
\]
We keep each interval \(I_\ell\) independently with probability \(p=p_\nu\), and
define \(F\) to be the union of those \(J_\ell\) for which \(I_\ell\) was kept. Then
\[
        \mathbb E\mu(F)=p(1-2\theta).
\]

Let \(L=L_\nu\ge m\) so large that
\[
        M\left(\frac{a_1}{R}\right)^{L+1}
        \le
        \frac{\theta}{4k}.
\]
For each
\[
        h(t)=\sum_{j=m}^{\infty}c_jt^j\in\mathcal V_{m,R,M,b,B},
\]
define the Taylor truncation
\[
        h_L(t):=\sum_{j=m}^{L}c_jt^j.
\]
Then, since $\sum_{j\geq m+1}|c_j|R^j\leq M$, for \(0\leq t\leq a_1\) we have
\begin{equation}\label{eq:tail-small}
        |h(t)-h_L(t)|
        \le
        \sum_{j=L+1}^{\infty}|c_j|t^j
        \le
        M\left(\frac{a_1}{R}\right)^{L+1}
        \le
        \frac{\theta}{4k}.
\end{equation}
Since
\[
        k\asymp_{b,m}(a_1\delta)^{-m},
\]
we may choose \(L\) so that, with
\[
        r:=L-m+1,
\]
one has
\begin{equation}\label{eq:visible-dim}
r\ll_{m,R,M,b,\theta} 1+\frac{|\log\delta|}{\log(R/a_1)}.    
\end{equation}
Indeed, after cancelling the common factor \(a_1^m\), it is enough to force
\[
        \left(\frac{a_1}{R}\right)^{L+1-m}
        \ll_{m,R,M,b,\theta}
        \delta^m.
\]

Now fix \(\gamma\in[0,1]\). Suppose that for some
\(h\in\mathcal V_{m,R,M,b,B}\) one has
\[
        \gamma+h(a_i)\in F\qquad(1\le i\le n).
\]
For each \(i\), let \(\ell_i=\ell_i(h)\) be the unique integer in
\(\{0,\ldots,k-1\}\) such that
\[
        \gamma+h(a_i)\in J_{\ell_i}.
\]
Since \(J_{\ell_i}\) has distance at least \(\theta/k\) from the endpoints of
\(I_{\ell_i}\), and since inequality~\eqref{eq:tail-small} holds, the point
\[
        \gamma+h_L(a_i)
\]
lies in the same partition interval \(I_{\ell_i}\). Thus the interval pattern
\[
        \pi_\gamma(h):=(\ell_1(h),\ldots,\ell_n(h))
\]
is one of the interval patterns arising from the truncated polynomials
\(h_L\).

We now count the possible interval patterns of the truncated polynomials. The
coefficient vector
\[
        (c_m,\ldots,c_L)
\]
has dimension
\[
        r=L-m+1.
\]
For fixed \(i\), the interval containing
\[
        \gamma+h_L(a_i)
\]
is determined by which partition endpoints lie to its left. Equivalently, it is
determined by the truth values of the affine inequalities
\[
        \gamma+c_ma_i^m+c_{m+1}a_i^{m+1}+\cdots+c_La_i^L
        \ge
        \frac{\ell}{k}.
\]
For fixed \(i\), as \(h\) ranges over \(\mathcal V_{m,R,M,b,B}\), the set of
possible values of \(h_L(a_i)\) is contained in an interval of length
\(O_{R,M}(a_i^m)\), hence in an interval of length \(O_{R,M}(a_1^m)\). Since the
partition endpoints are spaced by \(1/k\), only
\[
        O_{m,R,M,b}(ka_1^m+1)
        =
        O_{m,R,M,b}(\delta^{-m})
\]
values of \(\ell\) can affect the interval containing \(\gamma+h_L(a_i)\).
Therefore, as \(i\) ranges from \(1\) to \(n\), the possible interval patterns are
determined by at most
\[
        H\ll_{m,R,M,b}n\delta^{-m}
\]
affine inequalities in \(r\) variables.

Since \(R>0\) is fixed and \(a_1\to0\),
\[
        \log(R/a_1)=\log(1/a_1)+O_R(1)\sim \log(1/a_1).
\]
Thus the estimates~\eqref{eq:entropy-blocks} and~\eqref{eq:visible-dim} imply
\begin{equation}\label{eq:o(n)}
        r(\log n+|\log\delta|)=o(n),
\end{equation}
and in particular \(r=o(n)\). Let
\[
        H_0=C_{m,R,M,b}n\delta^{-m}
\]
be an upper bound for \(H\), with \(C_{m,R,M,b}\) chosen large enough that
\(H_0\ge n\). For all sufficiently large \(\nu\), we also have \(H_0\ge r\).
By Lemma~\ref{lem:threshold-count}, the number of possible interval patterns is
at most
\begin{equation}\label{eq:analytic-count}
 (r+1)H_0^r
        \le
        \exp\left(
        O_{m,R,M,b}\bigl(r(\log n+|\log\delta|)\bigr)
        \right).    
\end{equation}

By Claim~\ref{claim:analytic}, the values \(h(a_i)\) are separated by at least
\[
        \frac b2(a_1\delta)^m>\frac4k.
\]
Hence, whenever
\[
        \gamma+h(a_i)\in[0,1]\qquad(1\le i\le n),
\]
the intervals
\[
        I_{\ell_1(h)},\ldots,I_{\ell_n(h)}
\]
are distinct. Therefore, for any fixed interval pattern which can arise, the
probability that all corresponding intervals are kept is \(p^n\). It follows from
\eqref{eq:analytic-count} that, uniformly in \(\gamma\),
\[
        \mathbb P(\gamma\in X_\nu(F))
        \le
        \exp\left(
        O_{m,R,M,b}\bigl(r(\log n+|\log\delta|)\bigr)
        \right)p^n.
\]
After integrating over \(\gamma\in[0,1]\), the same bound holds for
\(\mathbb E\mu(X_\nu(F))\).

Put
\[
        S_\nu=C\,r(\log n+|\log\delta|),
\]
where \(C\) is larger than the implicit constant in the last display. Since
\(S_\nu=o(n)\) by estimate~\eqref{eq:o(n)}, choose
\[
        p_\nu=\exp\left(-\sqrt{\frac{S_\nu+1}{n}}\right).
\]
Then \(p_\nu\to1\), and
\[
        \exp(S_\nu)p_\nu^n
        =
        \exp\left(S_\nu-\sqrt{n(S_\nu+1)}\right)
        \to0.
\]
Thus
\[
        \mathbb E\mu(X_\nu(F))\to0.
\]
For all sufficiently large \(\nu\),
\[
        1-p_\nu<\varepsilon/4,
        \qquad
        \mathbb E\mu(X_\nu(F))<\eta/4.
\]
Since \(2\theta<\varepsilon/4\), we have
\[
        \mathbb E\mu([0,1]\setminus F)
        =
        1-p_\nu(1-2\theta)
        \le
        (1-p_\nu)+2\theta
        <
        \varepsilon/2.
\]
Markov's inequality gives
\[
        \mathbb P(\mu(F)\le 1-\varepsilon)<\frac12,
        \qquad
        \mathbb P(\mu(X_\nu(F))\ge\eta)<\frac14.
\]
Hence some realization satisfies
\[
        \mu(F)>1-\varepsilon,
        \qquad
        \mu(X_\nu(F))<\eta,
\]
as required.
\end{proof}

\begin{proof}[Proof of the analytic assertion in \cref{thm:randommain}]
We now prove the analytic assertion. Choose finite blocks
\(A_\nu=\{a_{\nu,1}>\cdots>a_{\nu,n_\nu}>0\}\) witnessing the additional hypothesis
\eqref{eq:entropy-main}; thus these blocks satisfy the hypotheses of
Lemma~\ref{lem:analytic-window}. Enumerate all analytic windows
\[
        \mathcal V_j=
        \mathcal V_{m_j,R_j,M_j,b_j,B_j},
\]
where
\[
        m_j\in\mathbb N,
        \qquad
        R_j,M_j,b_j,B_j\in\mathbb Q_{>0},
        \qquad
        b_j<B_j.
\]
This countable family covers every nonzero real-analytic function defined in a
neighbourhood of \(0\), after possibly restricting the function to a smaller
neighbourhood. Indeed, if
\[
        h(t)=c_mt^m+c_{m+1}t^{m+1}+\cdots
\]
is nonzero and real analytic near \(0\), then its power series has positive radius
of convergence. Choose a rational \(R>0\) smaller than this radius. Then
\[
        \sum_{j=m}^{\infty}|c_j|R^j<\infty.
\]
Choose rational numbers \(M,b,B>0\) such that
\[
        M>\sum_{j=m}^{\infty}|c_j|R^j,
        \qquad
        0<b<|c_m|<B.
\]
Then \(h\in\mathcal V_{m,R,M,b,B}\).

Choose positive numbers \(\varepsilon_j,\eta_j\) with
\[
        \sum_j\varepsilon_j<\frac14,
        \qquad
        \sum_j\eta_j<\frac18.
\]
Apply Lemma~\ref{lem:analytic-window} to each \(\mathcal V_j\). We obtain measurable sets
\(F_j\subset[0,1]\) and exceptional base-value sets \(X_j\subset[0,1]\) with
\[
        \mu(F_j)>1-\varepsilon_j,
        \qquad
        \mu(X_j)<\eta_j.
\]
Set again
\[
        F=\bigcap_{j=1}^\infty F_j,
        \qquad
        X=\bigcup_{j=1}^\infty X_j.
\]
Then \(\mu(F)>3/4\) and \(\mu(X)<1/8\). Choose an open set \(G\supset X\) with
\(\mu(G)<1/4\), and put \(E=F\setminus G\). Then \(E\) has positive measure.

Suppose, for contradiction, that \(f(A)\subset E\) for some nonconstant real-analytic
\(f\) defined near \(0\). Write
\[
        f(t)=\gamma+h(t),
        \qquad
        \gamma=f(0),
\]
where \(h\) is analytic and not identically zero. Since \(a_n\to0\) and
\(E\subset[0,1]\), we have \(f(a_n)\to f(0)=\gamma\), so \(\gamma\in[0,1]\). Choose
\(j\) such that \(h\in\mathcal V_j\). The block used in constructing \(F_j\) is a
subset of \(A\), and \(E\subset F_j\). Hence the definition of \(X_j\) gives
\(\gamma\in X_j\subset G\). Since \(G\) is open and \(f(a)\to\gamma\) along
\(a\in A\), \(a\to0\), we have \(f(a)\in G\) for all sufficiently small \(a\in A\),
contradicting \(f(A)\subset E=F\setminus G\). This proves the analytic assertion.
\end{proof}

\section{Proof of Corollary~\ref{cor:sufficient-criteria} and Theorem \ref{cor:lacunarySumsets}}
\subsection{Proof of Corollary~\ref{cor:sufficient-criteria}}
\begin{proof}[Proof of Corollary~\ref{cor:sufficient-criteria}]
We first prove that condition \(\mathrm{(i)}\) implies the analytic hypothesis in
\cref{thm:randommain}. Put
\[
        G_n=\max_{m\in T(n)}(a_m-a_{m+1}).
\]
Choose \(n_j\to\infty\) such that, with
\[
        a=a_{n_j},
        \qquad
        G=G_{n_j},
        \qquad
        \theta=G/a,
\]
one has \(\theta\to0\). Every interval contained in \([a/4,a]\) of length greater
than \(G\) contains at least one point of
\(\{a_{n_j},a_{n_j+1},\ldots,a_{\ell(n_j)}\}\); otherwise the sequence would jump
over such an interval, contradicting the definition of \(G\).

Choose disjoint intervals
\[
        J_r=\bigl(a-(4r+2)G,\,a-4rG\bigr],
        \qquad
        0\leq r<N,
\]
where
\[
        N=\left\lfloor\frac{a}{16G}\right\rfloor.
\]
Since \(\theta\to0\), we have \(N\to\infty\) and \(N\asymp \theta^{-1}\). For
\(j\) sufficiently large, all the intervals \(J_r\) are contained in \([a/2,a]\). Since
each \(J_r\) has length \(2G\), choose
\[
        y_r\in A\cap J_r\qquad(0\leq r<N).
\]
Then \(B_j=\{y_0>y_1>\cdots>y_{N-1}\}\subset A\), consecutive chosen points are
separated by at least \(2G\), and \(y_0\leq a\). Hence
\[
        \Delta(B_j)
        \geq
        \frac{2G}{a}=2\theta,
        \qquad
        |B_j|=N\asymp\frac1\theta.
\]
Therefore
\[
        |\log\Delta(B_j)|\ll\log\frac1\theta,
        \qquad
        \log N\ll\log\frac1\theta,
        \qquad
        N\gg\frac1\theta.
\]
It follows that \(-\log\Delta(B_j)=o(N)\). Moreover \(y_0\leq a_{n_j}\to0\), so
\(\log(1/y_0)\to\infty\), and for large \(j\),
\[
        \left(
        1+\frac{|\log\Delta(B_j)|}{\log(1/y_0)}
        \right)
        \bigl(\log N+|\log\Delta(B_j)|\bigr)
        \ll
        \left(\log\frac1\theta\right)^2
        =o\left(\frac1\theta\right)=o(N).
\]
Thus the analytic part of \cref{thm:randommain} applies.

It remains to show that \(\mathrm{(ii)}\) and \(\mathrm{(iii)}\) imply
\(\mathrm{(i)}\). If \(a_{n+1}/a_n\to1\), then for every \(\varepsilon>0\) and all
large \(m\),
\[
        a_m-a_{m+1}
        =a_m\left(1-\frac{a_{m+1}}{a_m}\right)
        \leq\varepsilon a_m.
\]
Hence for all large \(n\) and all \(m\geq n\),
\[
        a_m-a_{m+1}\leq\varepsilon a_m\leq\varepsilon a_n.
\]
Taking the maximum over \(m\in T(n)\) and then letting \(\varepsilon\downarrow0\)
gives \(\mathrm{(i)}\).

Finally assume \(\mathrm{(iii)}\), and write \(g_n=a_n-a_{n+1}\). Since
\[
        \frac{a_{n+1}}{a_n}=1-\frac{g_n}{a_n},
\]
there is a sequence \(n_j\to\infty\) such that \(g_{n_j}/a_{n_j}\to0\). Let
\(C_0=\max\{C,1\}\). For \(m\in T(n_j)\), we have \(m\geq n_j\), and therefore
\[
        g_m\leq C_0g_{n_j}.
\]
Consequently
\[
        \frac{\max_{m\in T(n_j)}(a_m-a_{m+1})}{a_{n_j}}
        \leq C_0\frac{g_{n_j}}{a_{n_j}}\to0.
\]
Thus \(\mathrm{(i)}\) holds, and the conclusion follows from the first part of the
proof.
\end{proof}

\subsection{Proof of Theorem~\ref{cor:lacunarySumsets}}

Let  $(b_i)$ and $(c_j)$ be two lacunary sequences with ratios $q_1,q_2$ and counting functions $I(W),J(W)$, respectively, and for $\sigma\in \{\pm 1\}$ define the (signed) near-additive energy \begin{align*}
    E_\sigma(W):=\#\Big\{(i,i',j,j'):i,i'\leq I(W),\ j,j'\leq J(W),\ |(b_i-b_{i'})+\sigma(d_j-d_{j'})|\leq \exp(-W)\Big\}.
\end{align*} The relevant estimate is as follows.

\begin{lemma}[\cite{cuellar2026nonuniversalitysumsetslacunarysequences}, Lemma 19] With $\max(q_1,q_2)<\frac{1}{2}$ there is a constant $C=C(q)$ such that for either $\sigma\in \{\pm1\}$ \begin{align*}
    E_\sigma(W)\leq CI(W)J(W)\quad \text{for all }W>0.
\end{align*}
    
\end{lemma}

One also defines the cross grid \begin{align*}
    G_\sigma(W):=\{b_i+\sigma c_j:i\leq I(W),j\leq J(W)\},\quad \sigma\in \{\pm 1\}.
\end{align*} The following proposition may then be deduced from their argument.

\begin{proposition}\label{prop:goodSeparation}
    Let $(b_i)$ and $(c_j)$ be two positive lacunary sequences with ratios $q_1,q_2$, respectively, and $\max(q_1,q_2)<\frac{1}{2}$. There are constants $c_0,c_1>0$, depending only on the lacunary ratios of the sequences, such that for any sufficiently large $W>0$ and either $\sigma\in \{\pm 1\}$, there is a $c_1\exp(-W)$-separated set \begin{align*}
        X_W\subset G_\sigma(W)
    \end{align*} with \begin{align*}
        |X_W|&\geq c_0I(W)J(W).
    \end{align*}
\end{proposition}

\begin{proof}
    Partition $[0,1]$ into $n:=\lceil \exp(W)\rceil$ intervals $I_1,\cdots, I_n$, each of length $\delta = \frac{1}{n}$. We have, by Cauchy-Schwarz, \begin{align*}
        I(W)J(W)=\sum_{k=1}^n\sum_{\substack{i\leq I(W)\\ j\leq J(W)}}\mathbf{1}(b_i+\sigma c_i\in I_k)\leq c(W)^{1/2}\Big(\sum_{k=1}^n\Big(\sum_{i,j}\mathbf{1}(b_i+\sigma c_i\in I_k)\Big)^2\Big)^{1/2},
    \end{align*} where $c(W)$ is the number of cells that are occupied; that is, \begin{align*}
        c(W):=\sum_{j=1}^n \mathbf{1}(\exists i,j:b_i+\sigma c_i\in I_j).
    \end{align*} Since \begin{align*}
        \sum_{k=1}^n\Big(\sum_{i,j}\mathbf{1}(b_i+\sigma c_i\in I_k)\Big)^2\leq E_\sigma(W)\leq CI(W)J(W)
    \end{align*} one then has that \begin{align*}
       \frac{I(W)J(W)}{C}\leq c(W).
    \end{align*} For each occupied cell $I_k$, choose one pair $(i,j)$ such that $b_i+\sigma c_i$ lands in $I_k$. The resulting set $(b_i+\sigma c_j)_{(i,j)\text{ chosen}}$ we call $X_W$.
\end{proof}

We are now able to prove Theorem \ref{cor:lacunarySumsets}.

\begin{proof}[Proof of Theorem \ref{cor:lacunarySumsets}] Let $(b_i)$ and $(c_j)$ be two positive lacunary sequences with ratios $q_1,q_2\in (0,1)$ and counting functions $I(W),J(W)$, respectively, such that \begin{align}\label{eq:growthCondition}
    \limsup_{W\rightarrow \infty}\frac{I(W)J(W)}{W}=\infty
\end{align} We first thin these sequences and pass to a tail, in a way that will not change (\ref{eq:growthCondition}). If $q_1\geq \frac{1}{2}$ then set $m:=\lceil \frac{\log q_1}{\log(1/4)}\rceil$ and pass to the subsequence $(b_{mi})_{i\in \mathbb{N}}$; similarly, if $q_2\geq \frac{1}{2}$, do likewise. This will not change (\ref{eq:growthCondition}), since the counting functions change by at most a constant factor. Passing to a fixed tail also will not change (\ref{eq:growthCondition}): if one starts at the $R$th and $S$th element of $(b_i)$ and $(c_j)$, respectively, then the new counting functions are $I(W)-R$ and $J(W)-S$. We claim that \begin{align*}
    \limsup_{W\rightarrow \infty}\frac{(I(W)-R)(J(W)-S)}{W}=\infty
\end{align*} follows from (\ref{eq:growthCondition}). Indeed, if $I(W),J(W)$ are both unbounded, then this is obviously true; if at least one is bounded, then since $\max(I(W),J(W))\leq C_{\max\{q_1,q_2\}}W$ for the lacunary sequences $(b_i)$ and $(c_j)$, one has a violation to (\ref{eq:growthCondition}).

With this in mind, we may assume that $\max(q_1,q_2)<\frac{1}{2}$. We now construct the blocks to apply Theorem \ref{thm:randommain}. Let $k\in \mathbb{N}$; since $(b_i),(c_j)$ tend to zero we may choose tail indices $R_k,S_k$ such that \begin{align*}
    b_{R_k}+c_{S_k}<2^{-k}.
\end{align*} Since \begin{align*}
    \limsup_{W\rightarrow \infty}\frac{(I(W)-R_k)(J(W)-S_k)}{W}=\infty
\end{align*} we may find some $W_k$ such that the following three points hold: \begin{enumerate}
    \item \begin{align*}
    \frac{(I(W_k)-R_k)(J(W_k)-S_k)}{W_k}\geq 8c_0^{-1}k,
\end{align*} where $c_0$ is the constant in Proposition \ref{prop:goodSeparation}

\item $W_k$ is large enough to appear in Proposition \ref{prop:goodSeparation}

\item $I(W_k)\geq 2R_k$, $J(W_k)\geq 2S_K$.
\end{enumerate} By Proposition \ref{prop:goodSeparation}, there are $c_0,c_1>0$ depending only on $q_1,q_2$, such that there is a $c_1\exp(-W_k)$-separated set \begin{align*}X_{W_k}\subset \{b_i+\sigma c_j:R_k\leq i\leq I(W_k),S_k\leq j\leq J(W_k)\}\end{align*} with $|X_{W_k}|\geq c_0I(W_k)J(W_k)$. Thus, \begin{align*}
    N_k:=|X_{W_k}|&\geq c_0I(W_k)J(W_k)\\ &= c_0\Big(1-\frac{R_k}{I(W_k)}\Big)^{-1}\Big(1-\frac{S_k}{J(W_k)}\Big)^{-1}(I(W_k)-R_k)(J(W_k)-S_k) \\ &\geq 32kW_k.
\end{align*} We have that \begin{align*}
    \Delta_k:=\min_{\substack{x,y\in X_{W_k}\\ x\neq y}}\frac{|x-y|}{\max X_{W_k}}\geq \exp(-W_k)
\end{align*} and so $-\log(\Delta_k)=o(N_k)$. Taking $k\rightarrow \infty$, we may then apply Theorem \ref{thm:randommain} to produce the desired result.
\end{proof}

\section{Degree-$d$ polynomial non-universality for $(d+3)A$}\label{sec:fixedDegreeNU}

Let us discuss the roadmap for the proof of Theorem~\ref{thm:sumsetNonuniversality}. Once and for all, choose distinct fixed points $b_1,\cdots,b_d\in A$. Since any polynomial $P$ of degree at most $d$ has at most $d-1$ critical points, any such polynomial cannot be flat on all of these fixed points. Restricting to a compact box of coordinates, one has uniformity in this observation, in that one of these points $b_\ell$ has a neighborhood on which $P$ roughly preserves distances.

Inside this neighborhood, we use points \begin{align*}
    b_\ell+a_{i_1}+\cdots+a_{i_{d+2}},\quad i_1<\cdots <i_{d+1}
\end{align*} with rapidly decreasing elements $a_j\in L\subset A$. For a fixed $b_\ell$, there are $\binom{|L|}{d+2}$ such points.

A degree-$d$ polynomial on $b_\ell+(d+2)L$ has $O(|L|^{d+1})$ degrees of freedom, and so an approximation net for the total space of configurations has size $O(|L|^{d+1}\log |L|)$, whereas the number of separated points is of order $|L|^{d+2}$. This gap between the two quantities allows one to construct a small random set that hits every configuration. Finally, a slicing argument gives the full set $E$.

\begin{definition}
    For a polynomial $P$ with degree at most $d$, we say that $b\in \{b_1,\cdots, b_d\}$ is a \textit{noncritical point} if $P'(b)\neq 0$.
\end{definition}

We first specialize to the compact family \begin{align*}
    \mathcal{P}:=\Big\{P(y)=\sum_{k=1}^d t_k y^k:1\leq \max_{1\leq k\leq d}|t_k|\leq 2\Big\}.
\end{align*}

\begin{lemma}\label{lem:goodSeparationNoncriticalPoints}
    There are constants $0<c\leq C$ and $0<\eta\leq 1$ such that every $P\in \mathcal{P}$ has a noncritical point $b_P=b(P)\in \{b_1,\cdots, b_d\}$ satisfying \begin{align*}
        c\leq |P'(b_P+s)|\leq C\quad (0\leq s\leq \eta).
    \end{align*} For any such point, \begin{align*}
        c|s-s'|\leq |P(b_P+s)-P(b_P+s')|\leq C|s-s'|\quad (s,s'\in [0,\eta]).
    \end{align*}
\end{lemma}

\begin{proof}
    The function $\phi:\mathcal{P}\rightarrow \mathbb{R}$ taking $P$ to $\max_\ell |P'(b_\ell)|$ is strictly positive: otherwise $P'$ would have $d$ distinct zeroes and degree at most $d-1$, forcing $P'=0$, a contradiction. By compactness, then, \begin{align}\label{eq:1.2.1}
        \min_{P\in \mathcal{P}}\max_\ell |P'(b_\ell)|=c'>0.
    \end{align} The quantities $|P'|$ and $|P''|$ are uniformly bounded on $[0,2]$ for $P\in \mathcal{P}$. Choose an upper bound $C\geq c'$ there. Uniform continuity of $P'$ gives $0<\eta\leq 1$ such that \begin{align}\label{eq:1.2.2}
        |P'(b_\ell+s)-P'(b_\ell)|\leq c'/2,\quad (P\in \mathcal{P},\ 1\leq \ell\leq d,\ 0\leq s\leq \eta)
    \end{align}

    Given $P\in \mathcal{P}$, we may find some $\ell$ such that $|P'(b_\ell)|\geq c'$ by (\ref{eq:1.2.1}), and so by (\ref{eq:1.2.2}), \begin{align*}
        |P'(b_\ell+s)|\geq |P'(b_\ell)|-|P'(b_\ell+s)-P'(b_\ell)|\geq c'/2,\quad (0\leq s\leq \eta).
    \end{align*} Setting $c=c'/2$ then provides the first part of the lemma. The mean value theorem then gives the second part.
\end{proof}

We now set up the family of points. Choose a lacunary family $L\subset A$ of points, disjoint from $\{b_1,\cdots, b_d\}$ in the following manner: \begin{align*}
    L=\{a_J<a_{J-1}<\cdots<a_1<\eta/2\},\quad 0<a_{j+1}\leq a_j/100\quad (1\leq j<J). 
\end{align*} For $I\subset [J]$ with $|I|=d+2$, set $s_I:=\sum_{j\in I}a_j$. Note that each $s_I$ lies in $[0,\eta]$. If $I\neq I'$, by considering their first differing index and the lacunary nature of $L$, the set $(s_I)_{I\subset [J]}$ has $\binom{J}{r}$ distinct elements, separated by at least $\Omega(a_J)$. For a sufficiently large constant $C_d\geq 1$ depending on $d$, we define the lacunary vector $v\in \mathbb{T}^J$ by \begin{align*}
    v:=\frac{1}{C_d}\begin{pmatrix}
        a_1 \\ a_2 \\ \vdots \\ a_J
    \end{pmatrix}.
\end{align*}

 For $P\in \mathcal{P}$, let $b_P\in  \{b_1,\cdots, b_d\}$ denote a choice of noncritical point given by Lemma \ref{lem:goodSeparationNoncriticalPoints}. Our total space of configurations will be \begin{align*}
    \mathcal{C}:=(z+(P(b_P+s_I)-P(b_P))v)_{\substack{I\subset [J]\\ |I|=d+2}},\quad z\in \mathbb{T}^J, P\in \mathcal{P}.
\end{align*} For each fixed $(z,P)$ pair, we view this as a list of vectors in $\mathbb{T}^J$. Equivalently, \begin{align*}
    \mathcal{C}=\bigcup_{\substack{z\in \mathbb{T}^J\\ P\in \mathcal{P}}}\mathcal{C}_{z,P},\quad\text{where}\quad \mathcal{C}_{z,P}:=(z+(P(b_P+s_I)-P(b_P))v)_{\substack{I\subset [J]\\ |I|=d+2}}.
\end{align*} We endow the space of configurations with an appropriate metric by setting \begin{align*}
    d(\mathcal{C}_{z,P},\mathcal{C}_{z',Q}):=\max_{\substack{I\subset [J]\\ |I|=d+2}}\Big\|(z+P(b_P+s_I)-P(b_P))v-(z'+Q(b_{Q}+s_{I})-Q(b_Q))v\Big\|_{L^\infty(\mathbb{T}^J)}.
\end{align*}

We first claim that $\mathcal{C}$ is adequate for the purposes of representing every degree-$d$ configuration.

\begin{lemma}\label{lem:sufficientCovering}
    Suppose there exists a set $H\subset \mathbb{T}^J$ such that, for each $z\in \mathbb{T}^J$ and $P\in \mathcal{P}$, there exists some $I\subset [J]$ with $|I|=d+2$ such that \begin{align*}
        z+(P(b_P+s_I)-P(b_P))v\in H.
    \end{align*} Then, for every $x\in \mathbb{T}^J$ and $P\in \mathcal{P}$ there exists some $y\in (d+3)A$ such that \begin{align*}
        x+P(y)v\in H.
    \end{align*}
\end{lemma}

\begin{proof}
    Given $x\in \mathbb{T}^J$ and $P\in \mathcal{P}$, set $z=x+P(b_P)v$. Then, by assumption there exists some $I\subset [J]$ with $|I|=d+2$ such that $x+P(b_P+s_I)v\in H$. Since $b_P+s_I\in (d+3)A$, setting $y=b_P+s_I$ provides the result.
\end{proof}

We now show we have the desired separation within each configuration.

\begin{lemma}\label{lem:wellSeparation}
    There exists an absolute constant $\rho>0$ depending only on $c$ and $C$ in Lemma \ref{lem:goodSeparationNoncriticalPoints} such that, for fixed $z\in \mathbb{T}^J$ and $P\in \mathcal{P}$, and letting $I$ vary over subsets of $[J]$ with size $d+2$, the points \begin{align*}
        z+P((b_P+s_I)-P(b_P))v
    \end{align*} are $\rho$-separated.
\end{lemma}

\begin{proof}
    Take $I\neq I'$ with $|I|=|I'|=d+2$. Let $m$ be their first differing index of $I$, $I'$; then by Lemma \ref{lem:goodSeparationNoncriticalPoints} and the definition above, \begin{align*}
        c|s_I-s_{I'}|\leq |P(b_P+s_I)-P(b_P+s_{I'})|\leq C|s_I-s_{I'}|.
    \end{align*} Since $|s_I-s_{I'}|=\Theta(a_m)$, by definition of $v$ one has that \begin{align*}
        \|P(b_P+s_I)-P(b_P+s_{I'})\|_\infty=\Omega_{c,C}(1).
    \end{align*} Set $\rho$ to be the lower bound here.
\end{proof}

We now want to control the space of configurations $\mathcal{C}$ with a finite number of representatives, with respect to the aforementioned metric. Let us define an $h$-net of $\mathcal{C}$ to be a finite subset $\mathcal{C}_0\subset \mathcal{C}$ of representatives, such that \begin{align*}
    \mathcal{C}_{z,P}\in \mathcal{C}\implies \text{ there exists some }\mathcal{C}_{z',Q}\in \mathcal{C}_0\text{ such that }d(\mathcal{C}_{z,P},\mathcal{C}_{z',Q})\leq h.
\end{align*}

\begin{lemma}\label{lem:hNet}
    For $h\in (0,1)$, the entire space of configurations $\mathcal{C}$ admits an $h$-net $\mathcal{C}_0$ of size \begin{align*}
        |\mathcal{C}_0|\ll_d h^{-O_d(J^{d+1})}.
    \end{align*}
\end{lemma}

\begin{proof}
    Fix some $b\in \{b_1,...,b_d\}$; by the union bound it suffices to show that the set \begin{align*}
        (z+(P(b+s_I)-P(b))v)_{\substack{I\subset [J]\\ |I|=d+2}},\quad z\in \mathbb{T}^J,P\in \mathcal{P}
    \end{align*} admits an $h$-net of the same size (up to a constant depending on $d$). We may expand \begin{align*}
        P(b+s_I)-P(b)=\sum_{k=1}^d c_k(b,P) (s_I)^k
    \end{align*} for some coefficients $c_1(b,P),...,c_k(b,P)$ depending only on $b$ and $P$, and so \begin{align*}
        P(b+s_I)-P(b)=\sum_{k=1}^dc_k(b,P)\Big(\sum_{j\in I}a_j\Big)^k=\sum_{k=1}^d \sum_{j_1,...,j_k\in I}c_k(b,P) a_{j_1}\cdots a_{j_k}.
    \end{align*}  Consequently, the $i$th entry of $z+(P(b+s_I-P(b))v$ is \begin{align}\label{eq:coordinateExpansion}
        z_i+\sum_{k=1}^d \sum_{j_1,...,j_k\in I}c_k'(b,P)\frac{a_{j_1}\cdots a_{j_k}}{a_i},\quad c_k':=\frac{c_k}{C_d}.
    \end{align} Consider the map \begin{align*}
        \phi:(P,z)\mapsto \Big(z_i, c_k'(b,P)\frac{a_{j_1}\cdots a_{j_k}}{a_i}\Big)_{\substack{1\leq i,j_1,...,j_k\leq J\\ 1\leq k\leq d}}
    \end{align*} (at possible abuse of notation, the image of $\phi$ here is a single vector). This image has $O(J^{d+2})$ coordinates, and so lies in $\mathbb{T}^M$ with $M=O(J^{d+2})$.

    For a constant $c_0(d)$ to be determined, partition this image space $\mathbb{T}^M$ into $\lceil \frac{c_0(d)}{h}\rfloor^M$ boxes by partitioning each coordinate into $\lceil \frac{c_0(d)}{h}\rceil$ boxes. For each box $B$ in the image space, choose some representative $(P,z)$ such that $\phi(P,z)\in B$. In such a way, every $(P,z)$ in the preimage has a representative $(Q,z')$ such that $\phi((P,z))$ and $\phi((Q,z'))$ lie in the same box. The map (\ref{eq:coordinateExpansion}) is linear in the entries of $\phi((P,z))$, and so by definition of the chosen representative one obtains the desired bound by applying the triangle inequality and choosing $\mathcal{C}_0$ to be the set of representatives.
\end{proof}

We can now construct an arbitrarily small set that hits every configuration.

\begin{proposition}\label{prop:hittingConfiguration}
For every $\delta\in (0,1)$, some $J\gg_d 1$, there exists a closed set $H\subset \mathbb{T}^n$ such that $\mu(H)\leq \delta$ and for every $z\in \mathbb{T}^J$, $P\in \mathcal{P}$, $z+P(b_P+s_I)v\in H$ for some $I\subset [J]$, $|I|=d+2$.
\end{proposition}

 \begin{proof}
     Let $\rho$ be the lower bound for the separation in Lemma \ref{lem:wellSeparation}, and fix an integer $k>1/\rho$. For a sufficiently large $J$ to be chosen, partition $\mathbb{T}^J$ into $k^J$ half-open cubes of side length $1/k$. Select each cube independently with probability $p$. Set $h:=\frac{1}{kJ}$ and let $H$ be the union of the closed $h$-neighborhoods of the selected cubes. This is a closed set, with expected measure \begin{align*}
         \mathbb{E}\mu(H)\leq pk^J(1/k+2h)^J=p(1+2hk)^J\leq ep.
     \end{align*} Take the net from Lemma \ref{lem:hNet} with radius $h$. By Lemma \ref{lem:wellSeparation}, each representative in this net has $\binom{J}{d+2}$ points in distinct cubes, since a cube has diameter at most $1/k<\rho$. The probability that a particular representative has no point in a selected cube is at most \begin{align*}
         (1-p)^{\binom{J}{d+2}}\leq \exp\Big(-p\binom{J}{d+2}\Big).
     \end{align*} The union bound then gives \begin{align*}
         \mathbb{P}(\text{some representative is missed})&\ll_d h^{-O_d(J^{d+1})}\exp\Big(-p\binom{J}{d+2}\Big) \\ & \ll_d \exp\Big(-p\binom{J}{d+2}+O_d(J^{d+1}\log J)\Big).
     \end{align*} Choosing $p=\delta/4e$ and taking $J$ sufficiently large, one may obtain $H$ such that $\mathbb{E}\mu(H)\leq \delta/2$ and \begin{align*}
         \mathbb{P}(\text{some representative is missed})\leq  1/8.
     \end{align*} By Markov's inequality, \begin{align*}
         \mathbb{P}(\mu(H)>\delta)\leq \frac{\mathbb{E} \mu(H)}{\delta}\leq \frac{1}{4}.
     \end{align*} Thus some choice of $H$ satisfies $\mu(H)\leq \delta$ and hits every representative. By construction of $H$, and the construction of the $h$-net, this set thus hits every configuration in $\mathcal{C}$.
     \end{proof}

     By Lemma \ref{lem:sufficientCovering}, this supplies that for every $x\in \mathbb{T}^J$ and $P\in \mathcal{P}$ there exists some $y\in (d+3)A$ such that $x+P(y)v\in H$. We may now construct a set for all degree-$d$ polynomial configurations.

\begin{proof}[Proof of Theorem \ref{thm:sumsetNonuniversality}]
    Take $\epsilon>0$. For each $n\in \mathbb{Z}$ we apply Proposition \ref{prop:hittingConfiguration} with $\delta_n= \frac{\epsilon}{12}2^{-|n|}$, to obtain a scale $J_n$, a set $H_n\subset \mathbb{T}^{J_n}$, and an associated vector $v_n$. We will use the normal product Haar probability measure on $G:=\prod_{n\in \mathbb{Z}}\mathbb{T}^{J_n}$. For $z=(z_n)_{n\in \mathbb{Z}}\in G$, define \begin{align*}
        E_z:= \{x\in [0,1]:z_n+x2^{-n}v_n\not\in H_n\text{ for every }n\in \mathbb{Z}\}. 
    \end{align*} Each $E_z$ is Borel. For fixed $x$ and $n$, letting $z$ vary uniformly over $G$, translation-invariance gives \begin{align*}
        \mathbb{P}(z_n+x2^{-n}v_n\in H_n)=\mu(H_n)\leq \delta_n.
    \end{align*} The union bound and Tonelli's theorem then produce that \begin{align*}
        \mathbb{E}\mu([0,1]\setminus E_z)\leq \sum_{n\in \mathbb{Z}}\delta_n= \frac{\epsilon}{4}.
    \end{align*} Consequently, some $z$ satisfies that $\mu(E_z)\geq 1-\epsilon$. Fix such a $z$ and write $E=E_z$. 

    Suppose a nonconstant polynomial $Q$ of degree at most $d$ satisfies $Q((d+3)A)\subset E$. By separating out its constant term and choosing a dyadic scale for its largest nonconstant coefficient, we may write \begin{align*}
        Q(y)=x_0+2^nP(y),\quad \text{for some }P\in \mathcal{P},n\in \mathbb{Z}.
    \end{align*} The definition of $E$ at index $n$ would then imply, for every $y\in (d+3)A$, \begin{align*}
        z_n+Q(y)2^{-n}v_n=z_n+x_02^{-n}v_n+P(y)v_n\not\in H_n.
    \end{align*} This is a violation of Proposition \ref{prop:hittingConfiguration}, and so by construction $Q((d+3)A)\not\subset E$.
\end{proof}

\section*{Acknowledgments}
The authors thank Xiang Gao and Joshua Zahl for helpful discussions. 

\bibliographystyle{abbrv}
\bibliography{main.bib}

\end{document}